\newenvironment{E}{\begin{equation}}{\end{equation}}
\def\proof{\noindent{\bf Proof: }}
\def\qed{ \hskip 20pt{\vrule height7pt width6pt depth0pt}\hfil}
\def\forb{{\hbox{forb}}}
\def\0{{\bf 0}}
\def\1{{\bf 1}}
\def\v{{\bf v}}

\documentclass[12pt]{article}
\renewcommand{\l}{\ell}

\def\Av{{\mathrm{Avoid}}}

\newtheorem{thm}{Theorem}[section]
\newtheorem{Lemma}[thm]{Lemma}

\newtheorem{cor}[thm]{Corollary}

\newcommand{\ncols}[1]{\| #1 \|}
\newcommand{\rf}[1]{(\ref{#1})}
\newcommand{\trf}[1]{Theorem~\ref{#1}}

\newcommand{\lrf}[1]{Lemma~\ref{#1}}

\newcommand{\srf}[1]{Section~\ref{#1}}

\usepackage{amssymb}
\usepackage{amsmath}
\title{Design Theory and some Forbidden Configurations}
\author{R.P.  Anstee\thanks{Research supported in part by
NSERC}, Farzin Barekat\thanks{Research supported in 2007 by NSERC of first author}, Zachary Pellegrin\thanks{Research supported in 2019 by NSERC USRA and NSERC of first author}
  \\Mathematics Department\\The University of British Columbia\\Vancouver,
B.C. Canada V6T 1Z2\\ {\small\texttt{anstee@math.ubc.ca}}
\\\mbox{\ }}

\begin{document}
\maketitle
\begin{abstract}
 In this paper we relate $t$-designs to a forbidden configuration problem in extremal set theory.  Let $\1_t\0_{\ell}$ denote a column of $t$ 1's on top of $\ell$ 0's.  Let $q\cdot \1_t\0_{\ell}$ denote the $(t+\ell)\times q$ matrix consisting of $t$ rows of $q$ 1's and $\ell$ rows of $q$ 0's. We consider  extremal problems for matrices avoiding certain submatrices. Let $A$ be a (0,1)-matrix  forbidding any $(t+\ell)\times(\lambda+2)$ submatrix $(\lambda+2)\cdot \1_t\0_{\ell}$. 
 Assume $A$ is $m$-rowed and only columns of sum $t+1,t+2,\ldots ,m-\ell$ are allowed to be repeated.   Assume that $A$ has the maximum number of columns subject to the given restrictions  assume $m$ is sufficiently large.  Then $A$ has each column of sum $0,1,\ldots ,t$ and $m-\ell+1,m-\ell+2,\ldots, m$ exactly once and, given the appropriate divisibility condition, the columns of sum $t+1$ correspond to a $t$-design with block size $t+1$ and parameter $\lambda$. The proof derives a basic upper bound on the number of columns of $A$ by a pigeonhole argument and then a careful argument, for large m,  reduces the bound by a substantial amount down to the value given by design based constructions. We  extend in a few directions.


\end{abstract}
\section{Introduction}

We explore a connection between  block designs and  extremal set theory. 
Combinatorial objects can be defined by forbidden substructures.  Let $[m]=\{1,2,\ldots ,m\}$. Also for any finite set $S$,   let $\binom{S}{k}$ denote all $k$-subsets of $S$. Thus $|\binom{[m]}{k}|=\binom{m}{k}$. An $m\times n$ (0,1)-matrix $A$ encodes a multiset ${\cal F}$ on $[m]$ consisting of $n$ sets (counted by multiplicity) where each column of $A$ is the incidence vector of a set in ${\cal F}$.   In this paper we consider certain (0,1)-submatrices (called \emph{configurations})  as the forbidden substructures of interest.    

 Consider a multiset ${\cal F}$ where each element of ${\cal F}$ is a called a \emph{block}. We say ${\cal F}$ is 
$t$-$(m,k,\lambda)$ \emph{design} if each set $B\in{\cal F}$ has
$B\in \binom{[m]}{k}$ (namely blocks of size $k$) and for each $t$-set $S\in\binom{[m]}{t}$ there are exactly $\lambda$ sets (blocks) in ${\cal F}$ containing $S$, with sets counted according to multiplicity .  It is usual in the study of block designs to use $v$ instead of  $m$ and to allow repeated blocks.  Recent results of Keevash \cite{K} yield that (for fixed $t,k,\lambda$ and $m$ large) $t$-designs exist assuming the easy \emph{divisibility conditions}: 
\begin{E}\binom{k-i}{t-i}\hbox{ divides } \lambda\binom{m-i}{t-i}\hbox{ for }i=1,2,\ldots , t-1.\label{easydivisibility}\end{E} 
Moreover, Keevash shows that we can require that there are no repeated blocks i.e. ${\cal F}$ is a set.   In that we case we call ${\cal F}$ a \emph{simple} design.

\begin{thm}(Keevash \cite{K}) Let $t,\lambda, k$ be given. Assume $m$ is sufficiently large and satisfies the divisibility conditions  \rf{easydivisibility}.  Then a simple 
$t$-$(m,k,\lambda)$ design exists. \label{keevash}\end{thm}

Our initial investigations \cite{AB} only considered $2$-$(m,3,\lambda)$ designs because we needed simple designs.  A result of Dehon \cite{D} establishes the existence of simple 2-$(m,3,\lambda)$ designs. The  result of Keevash \cite{K} above establishes the  existence of simple $t$-designs, for $t\ge 3$ and it suggested seeking greater generality than in \cite{AB}. 

In the paper we use a (0,1)-matrix interpretation of sets. Let $\ncols{A}$ denote the number of columns of $A$ (counted by multiplicity if that is relevant).  For an $m$-rowed matrix $A$ and a set $S\subseteq [m]$, let $A|_S$ denote the submatrix of $A$ formed by the rows $S$. Let $\1_t\0_{\ell}$ denote the $(t+\ell)\times 1$ vector of $t$ 1's on top of $\ell$ 0's. For a $s\times 1$ vector $\v$, Let $t\cdot \v$ denote the $s\times t$ matrix of $t$ copies of $\v$.

\begin{thm}\label{design}Let $t,\lambda,k,m$ be given. Let $A$ be a (0,1)-matrix with column sums $k$. Assume that for each each $t$-set $S\in\binom{[m]}{t}$ that $A|_S$ contains $\lambda\cdot\1_t$ and does not contain $(\lambda+1)\cdot\1_t$.
Then $\ncols{A}={\lambda\binom{m}{t}}/{\binom{k}{t}}$ and $A$ is the incidence matrix of a $t$-$(m,k,\lambda)$ design. \end{thm}

\trf{design} corresponds to the usual definition of a design. We could state a version of this by only requiring $A|_S$ contains 
$\lambda\cdot\1_t$ but also requiring $\ncols{A}={\lambda\binom{m}{t}}/{\binom{k}{t}}$.

Our motivation for studying these problems came from 
extremal set theory. 
An $m\times n$ (0,1)-matrix $A$ can be thought of a multiset of $n$
subsets of $[m]$. For an $m\times 1$ (0,1)-column $\alpha$, we define 
\begin{E}I(\alpha)=\{i\in [m]\,:\, \alpha\hbox{ has }1\hbox{ in row }i\}.\label{salpha}\end{E}
From this we define the natural multiset system $\cal A$ associated with the matrix $A$:
\begin{E}{\cal A}=\{I(\alpha_i)\,:\, \alpha_i\hbox{ is  column }i\hbox{ of }A \}.\label{calA}\end{E}
Similarly, if we are given a multiset system ${\cal A}$, we can form its \emph{incidence}  matrix $A$, as long as we don't care about column order. We define a {\it simple} matrix $A$ as a (0,1)-matrix with no repeated columns. In this case ${\cal A}$ yields a {\it set} system and it is in this setting that extremal set theory problems are typically  stated.

The property of forbidding a submatrix is usually extended to forbidding any row and column permutation of the submatrix.  Let $A$ and $F$ be (0,1)-matrices. We say that $A$ has $F$ as a \emph{configuration} and write $F\prec A$ if there is a submatrix of $A$ which is a row and column permutation of $F$. For the configuration $F=s\cdot\1_{t}\0_{\ell}$ only row permutations actually matter but  we are motivated by the study of forbidden configurations where row and column permutations matter \cite{survey}.

\begin{thm}\label{designconfig}Let $t,\lambda, k,m$ be given.  Let $A$ be an $m$-rowed  (0,1)-matrix with column sums $k$.  Assume that $(\lambda+1)\cdot \1_t\not\prec A$.  Then   
$$\ncols{A}\le   {\lambda\binom{m}{t}}/{\binom{k}{t}}. $$
Moreover in the case of equality the columns of columns in $A$ form the incidence matrix of a $t$-$(m,k,\lambda)$ design. \end{thm}
\proof Let $A$ be a $m\times {\lambda\binom{m}{t}}/{\binom{k}{t}}$ incidence matrix of a $t$-$(m,k,\lambda)$ design. Then $A$ provides a construction yielding the lower bound. The upper bound follows from a straightforward pigeonhole argument. Each column of sum $k$ has $\binom{k}{t}$ $t$-subsets of rows containing $t$ 1's. We can have at most $\lambda$ such columns for a given $t$-set in $A$.  
\qed

\medskip

%The case $k=t+1$ follows  from the more general result %$\forb(m,\lambda\cdot \1_t)$ \cite{AF}. 
Our extremal matrices, under a forbidden configuration restriction, yield a design.
The following result is the analogue of \trf{designconfig} where the hypotheses are  altered and weakened. The  theorem  relates  interesting forbidden  configurations (i.e. $(\lambda+2)\cdot\1_t\0_{\ell}$) to designs. 
Some repeated columns are allowed in the manner of Design Theory investigations. 
Our main result is:
\begin{thm}\label{genldesignconfig}Let $t,\ell,\lambda$ be given with $t>\ell$.  Let $A$ be an $m$-rowed  (0,1)-matrix with no repeated columns of sum $0,1,\ldots ,t$ nor column sums $m-\ell+1,m-\ell+2,\ldots ,m$. Assume that $(\lambda+2)\cdot \1_t\0_{\ell}\not\prec A$.  Then there exist an $M$ so that for $m\ge M$,  
\begin{E}\ncols{A}\le  \sum_{i=0}^{t-1}\binom{m}{i}+ \left(1+\frac{\lambda}{t+1}\right)\binom{m}{t} +\sum_{i=m-\ell+1}^{m}\binom{m}{i}.\label{1..10bd} \end{E}
Moreover in the case of equality the columns of column sum $t+1$ in $A$ form the incidence matrix of a $t$-$(m,t+1,\lambda)$ design and all columns of sum $1,2,\ldots ,t$ and $m-\ell+1,m-\ell+2,\ldots ,m$ are present. \qed\end{thm}
Of course, a similar result holds for $\ell>t$. We chose the multiplier for $\1_t\0_{\ell}$ to be $\lambda+2$ so that we would end up with a $t$-$(m,t+1,\lambda)$ design and  connect with Design Theory notation. We will prove this in \srf{generalcases}. General Lemmas are in \srf{basicLemmas}. 

This specializes to a forbidden configuration result. Define  
$$\Av(m,F)=\{A\,:\,A\hbox{ is an }m\hbox{-rowed matrix with }F\not\prec A\}.$$ 
 Then our extremal set theory problem is:
$$\forb(m,F)=\max_A\{\ncols{A}\,:\,A\in\Av(m,F), A\hbox{ is simple }\}.$$ 
These problems have been extensively investigated \cite{survey}. Exact results have been rare for non-simple configurations  $F$. One exception is $\forb(m,q\cdot K_t)=\forb(m,q\cdot \1_t)$ where the design constructions achieve equality  \cite{AF}. We would like to handle all cases $F=(\lambda+2)\cdot(\1_t\0_{\l}$)  but this paper only succeeds when $t>\ell$ (or, of course, $t<\ell$ by taking (0,1)-complements) or the special case $t=\ell=2$ (see \trf{1100}).   Define $K_k$ to be the $k\times 2^k$ incidence matrix of all subsets of $[k]$ and define $K_k^s$ to be the $k\times \binom{k}{s}$ incidence matrix of $\binom{[k]}{s}$.   
We specialize \trf{genldesignconfig} to simple matrices to obtain the following.

\begin{cor}\label{exactconfigresult} Let $t,\ell,\lambda$ be given with $t>\ell$. Then there exists an $M$ so that for $m>M$:
\begin{E}\forb(m,(\lambda+2)\cdot \1_t\0_{\ell}) \le  
\sum_{i=0}^{t-1}\binom{m}{i}+ \left(1+\frac{\lambda}{t+1}\right)\binom{m}{t} +\sum_{i=m-\ell+1}^m\binom{m}{i}.\label{1..10bdv2}\end{E}
Equality is only achieved when (\ref{easydivisibility}) is satisfied (with $k=t+1$) and for a matrix $[K_m^0K_m^1\cdots K_m^t\,T\,K_m^{m-\ell+1}\cdots K_m^m]$ where $T$ is the incidence matrix of a simple \hfill \break
$t$-$(m,t+1,\lambda)$ design.\end{cor}

\vskip 10pt
We have some alternate constructions for the case of equality for small $m$. For $m=\lambda+t+\ell$ we have 
$A=[K_m^0K_m^1K_m^2\cdot K_m^tK_m^{t+1}K_m^{m-\ell+1}K_m^{m-\ell+2}\cdots K_m^m]$  $\in\Av(m,(\lambda+2)\cdot\1_t\0_{\ell})$.  Choose some subset $S\in\binom{[m]}{t}$. We check that there are exactly $m-t$ columns in $K_m^{t+1}$ that are all 1's on rows $S$ and one further column of $K_m^t$ that is all 1's on rows $S$.  Moreover on the remaining $m-t=\lambda+\ell$ rows we have $(\lambda+1)\0_{\ell}\not\prec K_{\lambda+\ell}^1$.     Thus our construction is in $\Av(m,(\lambda+2)\cdot\1_t\0_{\ell})$.  We check that the total number of columns is bigger than our design construction bound \eqref{1..10bd} by $\binom{m}{t+1}-\frac{\lambda}{t+1}\binom{m}{t}=\frac{\ell}{t+1}\binom{\lambda+t+\ell}{t}$
(given $m=\lambda+t+\ell$).    Thus we need some condition on $m$ being large in order to obtain our result.  In essence, the pigeonhole argument explored in
    Lemma~\ref{pigeon1..10} is insufficient to prove our bounds.
The following would be a version of \trf{genldesignconfig} more in keeping with \trf{design}.

\begin{thm}\label{design1..10} Let $t,\ell,\lambda$ be given with $t>\ell$.  There exists an $M$ so that for $m>M$, if $A$ is an $m\times n$ (0,1)-matrix with column sums in $\{t+1,t+2,\ldots ,m-1\}$ and  $A\in\Av(m,(\lambda+1)\cdot\1_t\0_{\ell})$ then
\begin{E}n\le \frac{\lambda}{t+1}\binom{m}{t}\end{E}
and we have equality if and only if the columns of $A$ correspond to the $(t+1)$-sets of a
$t$-$(m,t+1,\lambda)$ design. \qed\end{thm}

 The cases with $t=\ell$ would be more difficult since $\1_t\0_{t}$ is self-complementary (under (0,1)-complements) and matrices in $\Av(m,\1_t\0_{t})$ could easily have large column sums.   \trf{nonsimpleq1100} in \srf{extrastuff} is an example of this. For small $m$ such as $m=\lambda+3$, the construction
 $A=[K_m^0K_m^1K_m^2K_m^3K_m^{m-2}K_m^{m-1}K_m^m]$ 
$\in\Av(m,(\lambda+2)\cdot\1_2\0_2)$ and exceeds the bound (\ref{simpleq1100bd}) below much as is true above for $(\lambda+2)\cdot\1_t\0_{\ell}$. Since we were unable to generalize \trf{nonsimpleq1100} to $\1_t\0_t$ for $t\ge 3$, we will not prove \trf{nonsimpleq1100} here but state it for completeness in \srf{extrastuff} with the unrefereed proof in the  arXiv  \cite{AB}.

The proof of Theorem~\ref{nonsimpleq1100} uses Tur\'an's bound (\trf{turan}). If there were improvements on Tur\'an's bound \cite{C} for $t\ge 3$, they might help handle the configurations $(\lambda+3)\cdot \1_t\0_t$ but there has been no recent improvements in the bounds.

\section{Basic Lemmas for $\Av(m,(\lambda+2)\cdot\1_t\0_{\ell})$}\label{basicLemmas}
First we state an important result we use. 
The following bounds were proven by Tur\'an \cite{T} for $t=2$ and by de Caen \cite{C} for general $t$. Perhaps better bounds are possible for general $t$.

\begin{thm}\label{turan}Tur\'an Bounds. Let $k,t,m$ be given with $k\le m$. Let $G$ be a collection of $n$  distinct sets in $\binom{[m]}{t}$.  For $n$ satisfying
$$n\ge \binom{m}{t}-\frac{m-k+1}{m-t+1}\binom{m}{t}/ \binom{k-1}{t-1},$$ 
there exists a set $S\subset[m]$ with $|S|=k$ so that all $\binom{k}{t}$ $t$-subsets $\binom{S}{t}$ are in $G$. \qed\end{thm}

We follow some of the  arguments noted in \cite{AB} as well as new arguments to obtain
\trf{genldesignconfig}.  The Lemmas below  consider an $m$-rowed matrix $A\in \Av(m,(\lambda+2)\cdot\1_t\0_{\ell})$. 
%We will only be using these Lemmas when $\ell=1$ but felt that the extra generality is almost free and may be helpful in further investigations.  
Assume $A$ has the property that the column sums are restricted to $\{t,t+1,\ldots ,m-{\ell}\}$ and that columns of sum $t$ are not repeated.  
Note that columns of sum $1,2,\ldots ,t-1$ and $m-\ell+1,m-\ell+2,\ldots ,m$ do not contribute to the forbidden configuration 
$(\lambda+2)\cdot\1_t\0_{\ell}$.
Also note that if we allowed repeated columns of sum $t$, then we would get a less interesting result.
For $i=t,t+1$, let $a_i$ denote the number of columns of column sum $i$ in $A$ and let $a_{\geq t+2}$ denote the number of columns of column sum at least $t+2$ in $A$.  Given the nature of the extremal matrices yielding equality in the bound (\ref{1..10bd}), we expect $a_t=\binom{m}{t}$ and $a_{\ge t+2}=0$.
We will be proving \trf{genldesignconfig} by contradiction and will assume that 
\begin{E}a_t+a_{t+1}+a_{{\geq} t+2}> \left(1+\frac{\lambda}{t+1}\right)\binom{m}{t}.\label{atat+1at+21..10}\end{E}
Our first Lemma extends a pigeonhole principle.
\begin{Lemma}\label{pigeon1..10} Let $m,t,\ell,\lambda$ be given with $t>\ell$. Let $A$ be an $m\times n$ matrix with no $(\lambda+2)\cdot (\1_t\0_{\ell})$, columns sums in $\{t,t+1,\ldots ,m-\ell\}$ and with no repeated columns of sum $t$. Assume $m\ge t+\ell+\lambda+2$ and (\ref{atat+1at+21..10}).  Then
$$\binom{t}{t}\binom{m-t}{\ell}a_t+
\binom{t+1}{t}\binom{m-t-1}{\ell}a_{t+1}
+\binom{t+2}{t}\binom{m-t-2}{\ell}a_{\geq t+2}$$
\begin{E}\le \binom{m}{t+\ell}\binom{t+\ell}{\ell}(\lambda+1).\label{pigeonholeq1..10}\end{E}
There exists positive constants $c_1,c_2$ so that
\begin{E}\binom{m}{t}-c_1m^{t-1}\le a_t\le \binom{m}{t} \label{atineq1..10}\end{E}
\begin{E}\hbox{ and }\quad a_{\geq t+2}\le c_2m^{t-1}.\label{at+2ineq1..10}\end{E}\end{Lemma}
\proof We note that a column of column sum $k$ has $\binom{k}{t}\binom{m-k}{\ell}$ configurations $\1_t\0_{\ell}$ and note that $\binom{k}{t}\binom{m-k}{\ell}\ge
\binom{t+2}{t}\binom{m-t-2}{\ell}$ for $t+2\le k\le m-1$.  Counting the configurations $\1_t\0_{\ell}$ (which can appear on $t+\ell$ rows in up to $\binom{t+\ell}{\ell}$ orderings) and using the pigeonhole argument yields (\ref{pigeonholeq1..10})

For $m\ge 3\ell+t+1$ we have $\binom{t+1}{t}\binom{m-t-1}{\ell}<
\binom{t+2}{t}\binom{m-t-2}{\ell}$.
Hence
$$\binom{m-t}{\ell}a_t+(t+1)\binom{m-t-1}{\ell}(a_{t+1}+a_{\geq t+2})\le \binom{m}{t+\ell}\binom{t+\ell}{\ell}(\lambda+1).$$
From (\ref{atat+1at+21..10}), we have $a_{t+1}+a_{\geq t+2}\ge \left(1+\frac{\lambda}{t+1}\right)\binom{m}{t}-a_t$.
We substitute and obtain
$$\binom{m-t-1}{\ell}\binom{m}{t}(t+1+\lambda)-\binom{m}{t+\ell}\binom{t+\ell}{\ell}(\lambda+1)\hfil$$
$$\hfil\le 
\biggl((t+1)\binom{m-t-1}{\ell}-\binom{m-t}{\ell}\biggr)a_t,$$
We deduce that there is a constant $c_1$ (will depend on $\lambda,t,\ell$) so that first half of (\ref{atineq1..10}) holds. The second half of (\ref{atineq1..10}) follows from the fact that no column of sum $t$ is repeated.
In a similar way we have 
$$\binom{m-t}{\ell}a_t+(t+1)\binom{m-t-1}{\ell}\left(1+\frac{\lambda}{t+1}\binom{m}{t}-a_t-a_{\geq t+2}\right)$$
$$+\binom{t+2}{t}\binom{m-t-2}{\ell}a_{\geq t+2}, \le\binom{m}{t+\ell}\binom{t+\ell}{\ell}(\lambda+1)$$
and when we substitute the upper bound of (\ref{atineq1..10}), we  deduce that there is a constant $c_2$
(will depend on $\lambda,t,\ell$) so that (\ref{at+2ineq1..10}) holds.\qed

\medskip

Partition $A$ into three parts: $A_{t}$ consists of the columns of column sum $t$, $A_{t+1}$ is the  columns of column sum $t+1$ and $A_{\geq t+2}$ is the columns of column sum greater or equal than $t+2$. We  construct  ${\cal A}_t$, ${\cal A}_{t+1}$ from $A_t$ and $A_{t+1}$ using the notations of (\ref{salpha}) and (\ref{calA}). Note that ${\cal A}_t$ is a set given that there are no repeated columns of sum $t$ while ${\cal A}_{t+1}$ is a multiset.   Let $S=\{i_1,...,i_t\}\in\binom{[m]}{t}$. Then define:
\begin{E}\mu(S)=\left\{ \begin{array}{ll} 1 & \hbox{ if } S \in {\cal A}_{t}\\  0 &\hbox{ if } S \notin {\cal A}_{t} \\ \end{array} \right.\quad ,\qquad
E=\{S\in\binom{[m]}{t}\,:\,\mu(S)=0\} = \binom{[m]}{t}\setminus {\cal A}_t
\quad .\label{epsilon}\end{E}
Thus $E$ denotes the $t$-sets missing from ${\cal A}_t$. We expect $E=\emptyset$. Now 
\begin{E}a_t=\sum_{S\in\binom{[m]}{t}}\mu(S)=\binom{m}{t}-|E|.\label{a2eq110}\end{E}
We deduce from (\ref{atineq1..10}) that $|E|\le c_1m^{t-1}$. 
We  use hypergraph degree definitions applied to the multiset ${\cal A}_{t+1}$.
For $S\in\binom{[m]}{t}$, define 
\begin{E}d(S)=|\{x\in [m]\,:\,S\cup x\in{\cal A}_{t+1}\}|,\label{d}\end{E} 
where we count the sets $S\cup x$  with their multiplicity in ${\cal A}_{t+1}\}$. For example, if we have $t=2$ and sets $\{1,2,3\}$, $\{1,2,4\}$, $\{1,2,4\}$ in  ${\cal A}_{3}$ (${\cal A}_{3}$ may have repeated columns) then $d(\{1,2\})=3$. Thus
\begin{E}(t+1)\cdot a_{t+1}=\sum_{S\in\binom{[m]}{t}}d(S).\label{at+1eq}\end{E}
%Let 
%$${\cal U}(S)=\{r\in[m]\,:\, S\cup r \in {\cal A}_{t+1}\}.$$
 Since $m>\lambda+\ell+t+1$ and we are avoiding $(\lambda+2)\cdot (\1_t\0_{\ell})$ in $A_{t+1}$ then 
 $|\{ G \,:\,G\in {\cal A}_{t+1}\hbox{ and }S\subset G\} |<\lambda+\ell+1$,where we count $G$'s according to the multiplicity in ${\cal A}_{t+1}\}$.

\vskip 5pt
\begin{Lemma} Let $A$ satisfy hypotheses of \lrf{pigeon1..10}. Then 
\begin{E}d(S) \leq (\lambda+1)-\mu(S). \label{dijineq}\end{E}
\label{dij} \end{Lemma}

\proof Recall (\ref{d}) for which we are counting by multiplicity the $(t+1)$-sets containing a given $t$-set $S$. We proceed to a contradiction by assuming the opposite of (\ref{dijineq}), namely we have an $S\in\binom{[m]}{t}$ with $d(S)+\mu(S)\ge \lambda+2$.
Let $B_1,B_2,\ldots ,B_{\lambda+2}$ denote $\lambda+2$ sets in ${\cal A}_t\cup{\cal A}_{t+1}$, each containing the $t$-set $S$.  Thus there are at most $\lambda+2$ elements in $\cup_iB_i$ which are not already in $S$. Thus for $m\ge\lambda+2+t+\ell$, we will have $\ell$ elements of $[m]$ not in any $B_i$ yielding the configuration $(\lambda+1)\cdot\1_t\0_{\ell}$, a contradiction.  \qed

\medskip

Our hypothetical extremal construction for $A$ avoiding $(\lambda+2)\cdot \1_t\0_{\ell}$ is $[K_m^t|T]$ where $T$ is the incidence matrix of a $t$-$(m,t+1,\lambda)$ design. In that case 
 all $t$-sets are present exactly once and for any $t$-set $S$, the number of sets, apart from the $t$-set $S\in\binom{[m]}{t}$,  containing $S$ is $\lambda$ and they are of size $t+1$.  Let $Y$ denote these `typical' $t$-sets : 
\begin{E}Y=\{S\in\binom{[m]}{t}\,:\,d(S)=\lambda\hbox{ and }\mu(S)=1\}\label{Y}\end{E}
We wish to show $Y=\binom{[m]}{t}$ in our proof of \trf{genldesignconfig}. The following Lemma is a step in that direction.

\medskip

\begin{Lemma}\label{Y1..10}Let $A$ satisfy hypotheses of \lrf{pigeon1..10}. There exists a constant $c_3$ so that 
\begin{E}|Y|\ge \binom{m}{t}-c_3m^{t-1}\label{yineq}\end{E}\end{Lemma}
\proof
We  partition  the $\binom{m}{t}$ $t$-sets $S$ into 3 parts: $Y$, $E$ and the rest. By Lemma~\ref{dij}, for each $S\in E$ we have $d(S)\le \lambda+1$. Note that for $S\notin Y\cup E$, we have $\mu(S)=1$ and so $d(S)\le \lambda-1$ (else $S\in Y$). Thus from (\ref{at+1eq}) and \lrf{dij}, 
$$(t+1)a_{t+1}=\sum_{S\in\binom{[m]}{t}}d(S)\le\left(\lambda|Y|
+(\lambda+1)|E|
+(\lambda-1)\biggl(\binom{m}{t}-|Y|-|E|\biggr)\right)$$
Hence 
\begin{E}
a_{t+1}\le\frac{1}{t+1}\left(\lambda\binom{m}{t}+|E|-\binom{m}{t}+|Y|+|E|\right)
\label{a3ineq110extend}\end{E}
Substituting estimates of $a_t$, $a_{t+1}$, $a_{\geq t+2}$ from (\ref{atineq1..10}), (\ref{at+1eq}), (\ref{at+2ineq1..10}) into (\ref{atat+1at+21..10}), we have
$$\binom{m}{t}-|E|+  \frac{1}{t+1}\left(\lambda\binom{m}{t}+2|E|-\binom{m}{t}+|Y|\right)
+c_2m^{t-1}
>\left(1+\frac{\lambda}{t+1}\right)\binom{m}{t}$$
We deduce $\left(\frac{2}{t+1}-1\right)|E|+\frac{1}{t+1}|Y| +c_2m^{t-1}>\frac{1}{t+1}\binom{m}{t}$ and so there exists a constant $c_3=(t+1)c_2$ so that (\ref{yineq}) holds.\qed

\medskip

\begin{Lemma} Let $k$ be given.  Use the notations of \srf{basicLemmas}. 
There exists an $M$ so that for $m\ge M$, there exists a set of $k$ rows $B$ such that for any $t$-set $S\in\binom{B}{t}$ then $S\in Y$.
\label{consturan}\end{Lemma}
\proof
Form a $t$-hypergraph $G$ of $m$ vertices corresponding to the rows of $A$ and 
with edge $S$ if and only if $S\in Y$. Thus by Lemma~\ref{Y1..10}, the number of edges ($t$-sets)
of $G$ is at least $\binom{m}{t}-c_3m^{t-1}$. We apply \trf{turan}, by a result of de Caen \cite{C}. Thus there exists an $M$ so that for $m\ge M$, there is a $B\subset [m]$ with $|B|=k$  so for any $S\in \binom{B}{t}$ we have $S\in Y$. Hence for  $S\in \binom{B}{t}$ we have   $d(S)=\lambda$ and $\mu(S)=1$. \qed

\section{Exact Bound for $(\lambda+2)\cdot (\1_t\0_{\ell})$}\label{generalcases}

The following two lemmas provide useful counting inequalities.  Our main idea  is that if we have column sums $t$, whether repeated or not, that avoid $(\lambda+2)\cdot \1_{t-1}$ we may use a straightforward pigeonhole bound that the number of columns is at most a constant times $m^{t-1}$. 
\begin{Lemma} Let $A_{t+1}$ be the columns of sum $t+1$ in $A$. Given any row $r\in[m]$, let $A_{t+1}^r$ be the submatrix of $A_{t+1}$ formed by the columns having a $1$ in row $r$, then
$$a_{t+1}^r = ||A_{t+1}^r|| \leq \frac{\lambda+1}{t}\binom{m-1}{t-1}.$$
\label{at1rbound}\end{Lemma}
\proof
Since $m>t+l+\lambda+2$, any matrix with column sums $t+1$ containing $(\lambda+2)\cdot \1_t$ must also contain $(\lambda+2)\cdot \1_t\0_{\ell}$, therefore $A_{t+1}^r$ must avoid $(\lambda+2)\cdot \1_t$.

Since each column in $A_{t+1}^r$ has a $1$ in row $r$, on the remaining $m-1$ rows the matrix must avoid $(\lambda+2)\cdot \1_{t-1}$ since adding in row $r$ would create $(\lambda+2)\cdot \1_t$ and therefore $(\lambda+2)\cdot \1_t\0_{\ell}$.

The bound for a matrix of column sum $t$ avoiding $(\lambda+2)\cdot \1_{t-1}$ on $m-1$ rows is $\frac{\lambda+1}{t}\binom{m-1}{t-1}$ (see \trf{design1..10}) , thus
$$||A_{t+1}^r|| \leq \frac{\lambda+1}{t}\binom{m-1}{t-1} \qed$$

\medskip

\begin{Lemma} Let $A_{t+1}$ be the columns of sum $t+1$ in $A$.  Given any set of rows $R\in[m]$, let $A_{t+1}^R$ be the submatrix of $A_{t+1}$ formed by the columns having a $1$ in any row $r\in R$ then
$$a_{t+1}^R = ||A_{t+1}^R|| \leq |R|\cdot\frac{\lambda+1}{t}\binom{m-1}{t-1}.$$
\label{at1Rbound}\end{Lemma}
\proof
By  definition, $A_{t+1}^R$ is all columns which are in $A_{t+1}^r$ for some $r\in R$, therefore
$$||A_{t+1}^R|| \leq \sum_r ||A_{t+1}^r|| \leq |R|\cdot \frac{\lambda+1}{t}\binom{m-1}{t-1} \qed$$

\begin{Lemma} Let $R\in[m]$ be any set of rows of constant size $|R|=\rho$. Define $A_{t+1}^R$ according to the notation of \lrf{at1Rbound} and construct ${\cal A}_{t+1}^R$ from $A_{t+1}^R$ using the notations of (\ref{salpha}) and (\ref{calA}). Let
$$W_R = \{S\in\binom{[m]}{t}\, : \, \exists\, x\in[m]\text{ s.t. } S\cup x\in {\cal A}_{t+1}^R\}.$$
Let
$$Z_R = Y\setminus W_R\,\,.$$
Then, for $m$ large enough, $|Z_R| > 0$.\label{Zsize}\end{Lemma}
\proof
Each $t+1$-set in ${\cal A}_{t+1}^R$  contributes $t+1$ $t$-sets to $W_R$, therefore
$$|W_R| \leq (t+1)\cdot|{\cal A}_{t+1}^R| \leq (t+1)\rho\cdot \frac{\lambda+1}{t}\binom{m-1}{t-1}.$$
Thus there will exist a constant $c_4$ with 
$$|Z_R| = |Y|-|W_R| \geq \binom{m}{t}-c_3m^{t-1} - (t+1)\rho\cdot \frac{\lambda+1}{t}\binom{m-1}{t-1} \geq \binom{m}{t}-c_4m^{t-1}.$$
For $m$ sufficiently large, $\binom{m}{t}>c_4m^{t-1}$, therefore $|Z_R| > 0$ $\qed$

\medskip
In the following lemma, we use the result $Z_R\ne\emptyset$.

\begin{Lemma} $A$ has no column with fewer than $\lambda+\ell$ $0$'s\label{enough0}\end{Lemma}
\proof We assume for contradiction that there exists some column $\alpha$ with fewer than $\lambda+\ell$ $0$'s. Let $R\subset [m]$ be the set of rows on which $\alpha$ is $0$. By assumption $|R| < (\lambda+\ell)$ but we also note $R>\ell$ by assumptions on $A$. Construct $Z_R$ according to \lrf{Zsize}. Then since $|R|$ is bounded by a constant, by \lrf{Zsize}, $|Z_R|>0$.

Thus there exists some $t$-set $S\in Z_R$ such that $d(S) = \lambda$, $\mu(S) = 1$ since $Z_R\subseteq Y$ and using \rf{Y}. Each of the $\lambda+1$ columns contributing to this count has $0$'s in all rows of $R$ since this count consists of a column of sum $t$ which is $0$ outside of $S$ and columns of sum $t+1$ which do not contribute to $W_R$ and are therefore not in $A_{t+1}^R$. Thus these columns must have no $1$ in $R$. Since the columns of sum $t+1$ have no $1$ in rows $R$ we must have that $S\cap R = \emptyset$ and therefore the column of sum $t$ is also $0$'s in rows $R$. Also, since $\alpha$ is $1$ outside of $R$, $\alpha$ has only $1$'s in rows $S$ and is also $0$'s in rows $R$.

Therefore we have $\lambda$ columns of sum $t+1$ with $1$'s in $S$, a column of sum $t$ with $1$'s in $S$ and a column $\alpha$ has $1$'s in rows $S$. Thus we have $(\lambda+2)\cdot \1_t$ in rows $S$ in these $\lambda+2$ columns. Additionally as argued above, each of these columns has $0$'s in rows $R$. Recalling that $|R|\ge \ell$, this creates the forbidden object $(\lambda+2)\cdot \1_t\0_{\ell}$ on these columns, a contradiction.

Thus no such column $\alpha$ can exist and so $A$ has no column with fewer than $\lambda+\ell$ $0$'s. $\qed$

\begin{Lemma} Let $A\in\Av(m,(\lambda+2)\cdot\1_t\0_{\ell})$ with column sums in
$\{t,t+1,\ldots, m-\ell\}$. Assume no column of sum $t$ is repeated.  Then
$$||A|| \leq \left(1+\frac{\lambda}{t+1}\right)\binom{m}{t}$$
\label{upperbound}
\end{Lemma}
\proof
Assume for contradiction that $A$ exceeds this bound then $A$ must exceed the bound for avoiding $(\lambda+2)\cdot 1_t$, so must contain this object on some set of $\lambda+2$ columns. Let $D$ be the matrix of these columns. Therefore $A$ contains the rows of $1$'s of the forbidden object and it remains to show that $A$ must also contain the rows of $0$'s in appropriately chosen columns. If there are at least $\ell$ rows of $0$'s in $D$ then $A$ contains the forbidden object. Otherwise $D$ must have at most $\ell-1$ rows of 0's and all other rows have at least one $1$.

By \lrf{consturan}, with $k=(\lambda+2)\cdot (t+1)+\ell$, there exists some clique of rows $B\subset [m]$ with $|B| = (\lambda+2)\cdot (t+1)+\ell$ for which any $t$-set is in $Y$.   Recall that  $D$ must have at most $\ell-1$ rows of 0's and all other rows have at least one $1$ and hence $D|_B$ has at least $(\lambda+2)\cdot (t+1)+1$ 1's. By pigeonhole principle,
there must be some column $\beta$ in $D$ with $t+2$ $1$'s in the rows of $B$. Take any $t$-set from these $t+2$ rows. This set of $t$ rows of $B$ must be in $Y$. Therefore, since $\alpha$ has column sum at least $t+2$, there exist $\lambda$ columns of sum $t+1$ with $1$'s in these rows and a column of sum $t$ with $1$'s in these rows. These along with $\beta$ create the rows of $1$'s in the forbidden object.

All other rows of the column of sum $t$ are $0$ and the columns of sum $t+1$ have at most $\lambda$ rows in which they are not all $0$. By \lrf{enough0} $\beta$ has at least $\lambda+\ell$ $0$'s. Therefore, there are at least $\ell$ rows in which all of these columns are $0$ creating the forbidden object, a contradiction.

Thus $A$ must satisfy this bound. $\qed$

\medskip

\noindent{\bf Proof of \trf{genldesignconfig}}: The upper bound follows from \lrf{upperbound}. We now consider the case of equality. We can repeat the previous arguments with equality and \lrf{consturan} and \lrf{enough0} with be true. Using the same arguments as \lrf{upperbound} we see that if $A$ contains $(\lambda+2)\cdot \1_t$ and equals the bound then $A$ must contain our forbidden object $(\lambda+2)\cdot \1_t\0_{\ell}$. Therefore a matrix $A$ achieving this bound avoids $(\lambda+2)\cdot \1_t\0_{\ell}$ if and only if $A$ avoids $(\lambda+2)\cdot \1_t$. Thus in the case of equality we must have that $a_t = \binom{m}{t}$, $a_{\geq t+2} = 0$, $a_{t+1}=\frac{\lambda}{t+1}\binom{m}{t}$ and the columns of sum $t+1$ correspond to a $t-(m,t+1,\lambda)$ design. $\qed$

\section{The cases $F=q\cdot \1_1\0_1$, $F=q\cdot \1_2\0_2$ and further problems}\label{extrastuff}

Use the notation $[A|B]$ to denote the matrix formed by concatenating  $A$ with $B$.   A related problem is attempting to compute
$\forb(m,[\1_{t+\ell}|\lambda\cdot\1_t\0_{\ell}])$.
We have shown 
$\forb(m,[\1_{t+1}|2\cdot\1_t\0_{1}])
=\forb(m,3\cdot\1_t)$ \cite{AKa}. Note that $3\cdot\1_t\prec [\1_{t+1}|2\cdot\1_t\0_{1}]$.   At this point we do not know but conjecture that $\forb(m,[\1_{t+1}|\lambda\cdot\1_t\0_{1}])
=\forb(m,(\lambda+1)\cdot\1_t)$ for large $m$. Note that $K_4\in\Av(m,[\1_3|3\cdot \1_2\0_1)$ but
$\ncols{A}=1+\frac{5}{3}\binom{m}{2}+\binom{m}{1}+\binom{m}{0}$. We are using our intuition about what drives the bound 
$\forb$. 

The following results appear in the unrefereed manuscript \cite{AB}. In \cite{AFS} we showed that 
$$\left\lfloor\frac{q+1}{2}m\right\rfloor+2\le 
\forb(m, q\cdot(\1_1\0_1))\le
\left\lfloor\frac{q+1}{2}m+\frac{(q-3)m}{2(m-2)}\right\rfloor+2$$
 where the upper bound obtained by a pigeonhole argument is achieved for $m=q-1$ by taking $A=[K_m^0K_m^1K_m^2K_m^{m-1}K_m^m]$.   For $m$ with $m\ge \max\{3q+2, 8q-19\}$, we are able to show that the lower bound is correct and  slice $\frac{(q-3)m}{2(m-2)}\approx \frac{q-3}{2}$ off a pigeonhole bound. It is likely that our bound is valid for smaller $m>q-1$. The case $q=4$, is Lemma 3.1 in \cite{AKa} and took a page to establish.
The unrefereed manuscript \cite{AB} contains the following.
\begin{thm}\label{q10} Let $q\ge 3$ be given. Then for $m\ge \max\{3q+2, 8q-19\}$, \begin{E}\forb(m, q\cdot \1_1\0_1)=\lfloor\frac{q+1}{2}m\rfloor +2.
\hskip 1in\qed\label{q10bd}\end{E}
\end{thm}

The lower bound is easy. For $m$ even or $q-3$ even, let $G$ be a (simple) graph on $m$ vertices for which all the degrees are $q-3$ and for $m,q-3$ odd let $G$ be a graph for which $m-1$ vertices have degree $q-3$ and one vertex has degree $q-4$. Such graphs are easy to construct.
Let $H$ be the vertex-edge incidence matrix associated with $G$, namely for each edge $e=(i,j)$ of $G$, we add a column to $H$ with 1's in rows $i,j$ and 0's in other rows. Thus $H$ is a simple $m$-rowed matrix  with 
$\lfloor\frac{(q-3)m}{2}\rfloor$ columns each of column sum 2.
The simple matrix $A=[K_m^0\,K_m^1\,H\,K_m^{m-1}\,K_m^m]$ has $\lfloor\frac{(q+1)m}{2}\rfloor+2$ columns  and no configuration
$q\cdot (\1_1\0_1)$ which establishes $\forb(m,q\cdot (\1_1\0_1))\ge \lfloor\frac{(q+1)m}{2}\rfloor+2$.

To prove Corollary~\ref{1100} and Theorem~\ref{design1100}, we would prove the following:

\begin{thm}\cite{AB}\label{nonsimpleq1100} Let $\lambda>0$, $m$,  be given. Let $A$ be an $m\times n$ (0,1)-matrix so that no column of sum 0,1,2, $m-2$, $m-1$ or $m$ is repeated. Assume $A\in\Av(m,(\lambda+3)\cdot (\1_2\0_2)$. Then there exists a constant $M$ so that for $m>M$, 
\begin{E}n\le
  2+2m+\left(2+\frac{\lambda}{3}\right)\binom{m}{2}
  \label{q1100bd}\end{E}
with equality for $m\equiv 1,3(\hbox{mod }6)$.
If $A$ is an $m\times\forb(m,(\lambda+3)\cdot (\1_2\0_2))$ matrix with $m>M$ and $m\equiv 1,3(\hbox{mod }6)$, then $A$ consists of all possible columns of sum 0, 1, 2, $m-2$, $m-1$ and $m$ once each and there are two positive integers $a,b$ satisfying $a+b=\lambda$ with the columns of column sum 3 correspond to a $2-(m,3,a)$ design and the columns of sum $m-3$  correspond to the complements in $[m]$ of the blocks of a $2-(m,3,b)$ and $A$ has no further columns. \qed\end{thm}

Specializing to simple matrices we obtain the following:

\begin{cor}\cite{AB}\label{1100} Let $\lambda>0$ be given. There exists a constant $M=M(q)$ so that for $m>M$, 
\begin{E}\hbox{forb}(m, (\lambda+3)\cdot (\1_2\0_2))
\le 2+2m+\left(2+\frac{\lambda}{3}\right)\binom{m}{2},\label{simpleq1100bd}\end{E}
We have equality in (\ref{simpleq1100bd}) for $m\equiv 1,3(\hbox{mod }6)$. 
If $A$ is an $m\times\forb(m,(\lambda+3)\cdot (\1_2\0_2))$ simple matrix with $m>M$ and $m\equiv 1,3(\hbox{mod }6)$,  then there exist positive integers $a,b$ with $a+b=\lambda$ so that $A$ consists of all possible columns of sum 0, 1, 2, $m-2$, $m-1$, $m$ and  with the columns of column sum 3 correspond to a $2-(m,3,a)$ design and the columns of sum $m-3$  correspond to the complements in $[m]$ of the blocks of a $2-(m,3,b)$ design and $A$ has no further columns.
\qed\end{cor}

A design theory version of this is: 

\begin{thm}\cite{AB}\label{design1100} Let $\lambda$ and $m$ be given integers. There exists an $M$ so that for $m>M$, if $A$ is an $m\times n$ (0,1)-matrix with column sums in $\{3,4,\ldots ,m-3\}$ and  $A\in\Av(m,(\lambda+1)\cdot\1_2\0_2)$ then
\begin{E}n\le \frac{\lambda}{3}\binom{m}{2}.
\label{design1100bd}\end{E}
We have equality in (\ref{design1100}) if and only if there are positive integers $a,b$ satisfying $a+b=\lambda$ and there are $\frac{a}{3}\binom{m}{2}$ columns of $A$ of column sum 3 corresponding to the blocks of a
$2-(m,3,a)$ design  and there are $\frac{b}{3}\binom{m}{2}$ columns of $A$ of column sum $m-3$ of $m-3$-sets whose complements (in $[m]$) corresponding to the blocks of a
$2-(m,3,b)$ design.\qed\end{thm}

\end{document}